\documentclass[12pt]{amsart}

\usepackage{common}
\usepackage{mathrsfs}

\title{Generic wide indiscernible sets}

\author{Alexander Usvyatsov}

 \address{Alexander Usvyatsov\\
 Mathematics Department\\
 Hebrew University of Jerusalem\\
 91904 Givat Ram, Israel\\
 }

\address{Alexander Usvyatsov\\
    Universidade de Lisboa \\
  Centro de Matem\'{a}tica e Aplica\cc{c}\~{o}es Fundamentais\\
  Av. Prof. Gama Pinto,2\\
  1649-003 Lisboa \\
  Portugal}

\begin{document}

\begin{abstract}
We prove that a wide Morley sequence in a wide generically stable type is isometric to the standard basis of an $\ell_p$ space for some $p$. 
\end{abstract}

\maketitle

\section{Introduction}

Wide types in the context of continuous model theory of Banach structures were introduced by Shelah and the author in \cite{ShUs}. A wide type is in a sense a ``correct'' notion of a ``large''  definable (or rather, type-definable) object in this context: it is a type that contains the unit sphere of an infinite-dimensional subspace of the monster model. Results in \cite{ShUs} indicate that wide types are useful and rich objects. It is also clear that if the type is assumed to be stable, techniques of geometric stability theory yield interesting consequences on the geometric structure of the set of its realisations. However, in \cite{ShUs} only minimal wide stable types are investigated. One of the main results there is that the Morley sequence in such a type (which turns out to be also the unique sequence of its wide extensions) is isometric to the standard basis of a Hilbert space $\ell_2$. This structural result was used in order to prove that any non-separable non-separably categorical Banach structure is prime over an underlying Hilbert space. 

In this  note we initiate the study of the geometric structure of generic sequences in wide types in a much more general context. Our main result is that a wide Morley sequence in a wide generically stable type is isometric to the standard basis of an $\ell_p$ space for some $1\le p < \infty$. Since a Morley sequence in a minimal wide stable type is automatically wide, one can consider this result as a more general analogue of the theorem from \cite{ShUs} mentioned above. 

\section{Preliminaries}

Let us fix a monster Banach structure $\fC$.  

The reader may assume for simplicity that $\fC$ is the monster model (a ``big'' saturated model) of a continuous first order theory $T$, so that  $T$ extends an (incomplete) theory of (real or complex) normed spaces (see Example 4.5 in \cite{BU} for the discussion of several approaches to constructing such a theory). In particular, $\fC$ expands an underlying Banach space (i.e., it is a Banach space with uniformly continuous extra-structure). 

However, let us  remark briefly that for the purposes of our paper, $\fC$ does not need to be saturated, and that our proofs do not require compactness. So a reader familiar with more general contexts of model theory of metric structures, e.g., monster metric spaces \cite{ShUs837} or homogeneous metric 
monsters \cite{ShUs928} (which essemtially means that $\fC$ is a metric structure, which is also $(D,\lam)$-homogeneous \cite{Sh:3} for a ``big'' $\lam$) will see that our proofs easily carry over to these frameworks. Of course we still assume that $\fC$ is a Banach structure, that is, a continuous expansion of an underlying Banach space. 

All sets discussed in the paper are assumed to be subsets of $\fC$, and all models -- complete elementary submodels of $\fC$. In the (most interesting) case of $\fC$ saturated, this just means that a model is a model of $T$ in the sense of continuous logic (which  implies by definition  that it is complete).  

\subsection{Basic notation}
We denote the space of types over a set $B$ by $\tS(B)$, where as $\cS(B)$ stands for the ``unit ball'' of $B$, that is, the collection of all elements of $B$ of norm 1 (naturally, we will only use this notation in case $B$ is a subspace of $\fC$). 

Notations are more or less standard. Models will normally be denoted by $M,N, \ldots$, sets –– by $A,B, C, \ldots$. In this paper, $B$ normally denotes a  complete subspace of $\fC$. We denote elements and finite tuples by letters $a,b,c, \ldots$, By $a_{<n}$ we mean the sequence of tuples $\lseq{a}{i}{n}$, similarly for $a_{<\omega}$, etc.

\bigskip

Given a sequence $I = \lseq{a}{n}{\om} $ of elements (not tuples) in a Banach space $B$, a \emph{normalized block} of $I$ is an element $b = b(n,k)$ of the form $\sum_{i=0}^{k-1} r_ia_{n+i}$, where $n,k<\om$, $r_ i$ and $\|b\| = 1$. By \emph{sequence of normalized blocks} of $I$ we mean a sequence \lseq{b}{m}{\om}

\subsection{Blocks and sequences of blocks} 
\label{sub:blocks_and_sequences_of_blocks}

Let us fix some notation and terminology.

Let $k<\om$ be a natural number, $\bar{r} = \seq{r_0,\ldots,r_{k-1}} \in \setR^k$ a sequence of real numbers. Denote by $\rho_{k,\bar{r}}$ or simply by $\rho_{\bar{r}}$ the following term:
\[
  \rho_{k,\bar{r}}(x_0,\ldots,x_{k-1}) = \rho_{\bar{r}}(x_0,\ldots,x_{k-1}) = \rho_{\bar{r}}(\x) = \sum_{i=1}^k r_ix_i
\]

Given a global $M$-invariant type $\pp$, we denote by $\rho_{k,\bar{r}}(\pp)$ the type of $\rho_{k,\r}(a_1,\ldots,a_k)$, where 
$\lseq{a}{i}{\om}$ is some (=any) sequence in $\pp$ over some (=any) model $N$ saturated over $M$. More precisely, we denote by 
$\rho_{k,\bar{r}}(\pp)$ the unique global $M$-invariant extension of this type. 

In other words, $\rho_{k,\bar{r}}(\pp)$ is the type of a certain block (determined by $\rho$) of a generic sequence of realizations of $\pp$; and it is also a global $M$-invariant type. 

Similarly, given an $M$-indiscernible sequence $I = \lseq{a}{i}{\om}$, denote by $\rho_{k,\bar{r}}(I)$ the sequence of appropriate ``blocks'' from $I$, that is, 

\[
  \rho_{k,\bar{r}}(I) = \lseq{b}{n}{\om}
\] 

where $b_n = \rho_{k,\bar{r}}(a_{kn},a_{kn+1},\ldots,a_{k(n+1)-1})$

We call the collection of all types of the form $\rho_{k,\bar{r}}(\pp)$ the \emph{linear span of $\pp$}, denoted by $\Span(\pp)$. If $\pp = \Av(\cU)$ for some ultrafilter $\cU$ on $M$, we also denote $\Span(\pp)$ by $\Span(\cU)$. It is easy to verify that all elements of $\Span(\cU)$ are themselves finitely satisfiable in $M$, hence are average types of ultrafilters on $M$. If $\qq = \Av(\cV)$ and $\qq = \rho_{k,\bar{r}}(\pp)$, we will also write $\cV = \rho_{k,\bar{r}}(\cU)$. 

The \emph{closure} in the logic topology of the linear span of $\pp$ is called the \emph{closed linear span of $\pp$}, or just the \emph{closed span of $\pp$}, and is denoted by $\overline{\Span(\pp)}$.


\subsection{Standard sequences and sequence spaces} 
\label{sub:standard_sequences_and_sequence_spaces}

Fix a fragment $\Delta$ of the language $L$ and an $L$-structure $B$. 

\begin{dfn}
   Let $I = \lseq{a}{i}{\om}$ be a sequence $\Delta$-based on $B$, and let $1\le p < \infty$, $\eps>0$. 
   \begin{enumerate}
     \item We say that $I$ is $(1+\eps)$-equivalent to the standard basis of $\ell_p$ if for every $k<\om$ and $\r \in \setR^k$, we have that $$\|\rho_{k,\bar{r}}(a_0, \ldots, a_{k-1}\| \sim_{(1+\eps)} \sqrt[\leftroot{-2}\uproot{15}p]{\sum_{i=0}^{k-1}|r_i|^p}$$

     where $t \sim_{(1+\eps)} s$ means $$(1-\eps)t \le s \le (1+\eps)t $$ 

     \item We say that $I$ is $(1+\eps)$-equivalent to the standard basis of $c_0$ if for every $k<\om$ and $\r \in \setR^k$, we have that $$\|\rho_{k,\bar{r}}(a_0, \ldots, a_{k-1}\| \sim_{(1+\eps)} \max\set{|r_i| \colon i<k}$$

     \item We say that $I$ is $(1+\eps)$-equivalent \emph{over $B$} to the standard basis of $\ell_p$ if for every $k<\om$, $\r \in \setR^{k+1}$ with $r_k = 1$, and $b \in B$, we have that $$\|\rho_{k+1,\bar{r}}(a_0, \ldots, a_{k-1},b\| \sim_{(1+\eps)} \|\sqrt[\leftroot{-2}\uproot{15}p]{\sum_{i=0}^{k-1}|r_i|^p a_0 + b}\|$$

     \item We say that $I$ is $(1+\eps)$-equivalent \emph{over $B$} to the standard basis of $c_0$ if for every $k<\om$, $\r \in \setR^{k+1}$ with $r_k = 1$, and $b \in B$, we have that $$\|\rho_{k+1,\bar{r}}(a_0, \ldots, a_{k-1},b\| \sim_{(1+\eps)}  \|\max\set{|r_i| \colon i<k} a_0 + b\|$$

   \end{enumerate}
 \end{dfn}

We will say ``isometric'' for ``1-equivalent''.

\medskip

\section{Generic wide sequences}


\begin{dfn}\label{dfn:widesets}
\begin{enumerate}
\item
  Let $M$ be a model, $p \in \tS(M)$ a wide type. Given an indiscernible sequence $I=\lseq{a}{i}{\om}$ in $p$, we call $I$ a \emph{wide} indiscernible sequence if $I\subseteq \cS(B)$, where and $B$ an infinite dimentional subspace of $\fC$ so that the unit sphere of $B$ is contained in the set of realizations of $p$, that is., $I \subseteq \cS(B)\subseteq p^{\fC}$. 
 \item
 	We say that $I$ as above is a \emph{generic wide} indiscernible sequence in $p$ if in addition $\tp(a_n/Ma_{<n})$ is finitely satisfiable in $M$ (that is, $I$ is a coheir sequence over $M$).
\end{enumerate}
\end{dfn}

Note that if $p$ is stable (or just generically stable), then $I$ is wide generic if and only if it is a wide \emph{Morley} sequence: that is $\tp(a_n/Ma_{<n})$ does not fork over $M$ for all $n<\om$. 

\smallskip

A typical example of a generic wide sequence is the unique sequence of wide extensions (equivalently, the unique Morley sequence) in a \emph{minimal} wide stable type. Such sequences were studied extensively in \cite{ShUs}, and it was shown that any such sequence is 1-equivalent to the standard basis of $\ell_2$.

%
%
%

\begin{obs}\label{obs:blocksequence}
  Let $I$ be a wide sequence in a wide type $p$ over a set $A$. Then every sequence of normalized blocks of $I$ 
  is a sequence of realizations of $p$.
\end{obs}

\begin{thm}\label{thm:main_theorem}
  Let $M$ is a Banach structure, $p \in \tS(B)$ a generically stable wide type, and $I=\lseq{a}{i}{\om}$ a wide generic indiscernible set in $p$. Then $I$ is isometric (1-equivalent) to the standard basis of a basic sequence space. 
\end{thm}
\begin{prf}


  Let $J = \lseq{b}{i}{\om}$ be an sequence of normalized blocks of $I$ 
 By a well known characterization of  basic sequence spaces (which essentially goes back to  the classical paper of Bohnenblust \cite{Boh}), it is enough to show that $\tp(I) = \tp(J)$ (quantifier free type is enough). In fact, we will show that $\tp(I/M) = \tp(J/M)$

  We shall prove by induction on $n<\om$ that $\tp(a_{<n}/M) = \tp(b_{<n}/M)$; this clearly suffices. 
  
  For $n=1$ there is nothing to prove. Assume that we know $\tp(a_{<n}) = \tp(b_{<n})$ for some $n\ge 1$, and let $\sigma \in \Aut(\fC/M)$ be such that $\sigma(b_{<n}) = a_{<n}$.

  Since $b_n$ is a normalized block of $I$, we have (by e.g. Observation \ref{obs:blocksequence}) that $b_n \models p$. 

  Recall that we denote by $\text{range}(b_n)$ the indices of $I$ that are involved in the definition of $b_n$. 

  Since $\min\text{range}(b_n)>\max\text{range}(b_{i})$ for all $i<n$, and $I$ is a Morley sequence over $M$, we also have 
  $b_n \ind_M b_{<n}$. 


  Similarly, $a_n \models p$ and $a_n \ind_M a_{<n}$. Combining this with the induction hypothesis, by generic stability (stationarity of non-forking extensions), since $M$ is a model, the independent pairs $a_{<n}a_n$ and $b_{<n}b_n$ also have the same type over $M$. This finishes the proof.

\end{prf}

\begin{rmk}
	Since $c_0$ is not generically stable (that is, a generically stable type can not have a Morley sequence isometric to $c_0$), Theorem \ref{thm:main_theorem} shows that any such sequence is in fact isometric to the standard basis of some $\ell_p$ for $1\le p < \infty$.
\end{rmk}

\bibliography{common.bib}
\bibliographystyle{alpha}
\end{document}